\documentclass{amsart}
\usepackage[dvips]{graphicx}
\usepackage{verbatim}
\usepackage{amscd}
\theoremstyle{plain}
  \newtheorem{thm}{Theorem}[section]
  \newtheorem{lem}[thm]{Lemma}
  \newtheorem{cor}[thm]{Corollary}

\theoremstyle{definition}

\theoremstyle{remark}
  \newtheorem{rem}[thm]{Remark}
  \newtheorem*{ack}{Acknowledgments}
\numberwithin{equation}{section}
\date{\today}
\newcommand{\Q}{{\mathbb{Q}}}
\newcommand{\cA}{{\mathcal{A}}}
\newcommand{\cG}{{\mathcal{G}}}
\newcommand{\cK}{{\mathcal{K}}}
\newcommand{\cL}{{\mathcal{L}}}
\newcommand{\cM}{{\mathcal{M}}}
\newcommand{\cS}{{\mathcal{S}}}
\newcommand{\size}{{\operatorname{size}}}
\newcommand{\e}{{\varepsilon}}
\newcommand{\fig }[1]
  {\raisebox{- 4.7mm}{\includegraphics[scale=0.2]{#1.ps}}}
\newcommand{\figm}[1]
  {\raisebox{- 9  mm}{\includegraphics[scale=0.4]{#1.ps}}}
\newcommand{\figl}[1]
  {\raisebox{-13.5mm}{\includegraphics[scale=0.5]{#1.ps}}}
\begin{document}
\title{Finite type invariants of knots via their Seifert matrices}
\author{Hitoshi Murakami}
\address{
  Department of Mathematics,
  School of Science and Engineering,
  Waseda University
  Ohkubo, Shinjuku-ku, Tokyo, 169-8555, Japan
  \\
  and
  Mittag-Leffler Institute,
  Aurav{\"a}gen 17,
  S-182 62, Djursholm,
  Sweden
}
\email{hitoshi@uguisu.co.jp}
\author{Tomotada Ohtsuki}
\address{
  Department of Mathematical and Computing Sciences,
  Tokyo Institute of Technology,
  Oh-okayama, Meguro-ku, Tokyo 152-8552, Japan
}
\email{tomotada@is.titech.ac.jp}
\thanks{
This research is supported in part by Grand-in-Aid for
Scientific Research, The Ministry of Education, Science,
Sports and Culture.
The first-named author is also supported in part by Waseda University
Grant for Special Research Projects (No. 98A-623).
}
\begin{abstract}
We define a filtration on the vector space
spanned by Seifert matrices of knots
related to Vassiliev's filtration on the space of knots.
Further we show that
the invariants of knots derived from the filtration
can be expressed by coefficients of the Alexander polynomial.
\end{abstract}
\maketitle
%%%%
The theory of finite type invariants (Vassiliev invariants) for knots
was first introduced by V.~Vassiliev \cite{Vassiliev:1990} and
reformulated by J.S.~Birman and X.S.~Lin \cite{Birman/Lin:INVEM93}.
M.~Kontsevich defined the universal Vassiliev invariant
\cite{Kontsevich:1993, BarNatan:TOPOL95} by using iterated integral.
The invariant takes values in the linear combinations of chord
diagrams and one can use it to construct an isomorphism from
the space of all the Vassiliev invariants of degree $d$ to the chord
diagrams with $d$ chords modulo diagrams with more chords.
\par
D.~Bar-Natan \cite{BarNatan:TOPOL95} extended the notion of chord diagrams
allowing trivalent vertices, which we call web diagrams in this paper.
He showed that the space of (the linear combinations of) chord
diagrams modulo the four-term relation coincides with the space of web
diagrams modulo the AS, IHX and STU relations.
So a main interest in the theory of Vassiliev invariants is the
study of web diagrams.
\par
In this paper we consider a classical knot invariant, the
(S-equivalent classes of) Seifert matrices related to Vassiliev
invariants.
Let $\cK$ be the vector space over $\Q$ spanned by knots and let
$$
  \cK\supset\cK_1\supset\cK_2\supset\cK_3\supset\cdots
$$
be Vassiliev's filtration of $\cK$.
A Vassiliev invariant of degree $d$ is defined to be a map
$\cK/\cK_{d+1}\to\Q$.
Further let $\cS$ be the vector space spanned by Seifert matrices.
There is a natural map $s:\cK\to\cS$ which takes a knot to its Seifert
matrix.
We consider a filtration
$$
  \cS\supset\cS_1\supset\cS_2\supset\cS_3\supset\cdots
$$
of $\cS$ induced from Vassiliev's filtration by $s$.
Our motivation is that which finite type invariant factors
$\cS/\cS_{d+1}$.
Since the Alexander-Conway polynomial can be defined by using
Seifert matrices, the finite type invariants coming from the
Alexander-Conway polynomial factor $\cS/\cS_{d+1}$.
The main result of this paper is that these are all that factor
it.
\begin{ack}
This work was done during the authors were visiting the Mittag-Leffler
Institute.
They would like to thank the Institute for its hospitality.
They also thank Andrew Kricker for his helpful comments about
web diagrams and claspers related to the Alexander-Conway polynomial.
\end{ack}
\section{Statement of the result}
Let $\cM$ be the set of integer matrices of even size such that
$M-{}^t M$ is unimodular.
We also include the $0\times0$ matrix in $\cM$.
Two matrices in $\cM$ are called S-equivalent if one
can be obtained from the other by using the following three types of
transformations.
\begin{gather*}
    M\Leftrightarrow PM'{}^{t}P,
  \\
    M\Leftrightarrow
    \begin{pmatrix}
      0 & 0 & O \\
      1 & x & O \\
      O & C & M
    \end{pmatrix},
  \\
  M\Leftrightarrow
    \begin{pmatrix}
      0 & 1 & O \\
      0 & x & R \\
      O & O & M
    \end{pmatrix}.
\end{gather*}
Here $P$ is a unimodular matrix, $O$ is a zero matrix of suitable
size and $C$ and $R$ are some column and row matrices respectively.
We denote by $[M]$ the S-equivalence class of $M\in\cM$.
\par
Let $\cS$ be the vector space over $\Q$ spanned by S-equivalence
classes of matrices in $\cM$.
For a matrix $M \in \cM$ and integers $i_1,i_2,\dots,i_d$
with $1\le i_k\le\size(M)$ ($\size(M)$ is the size of $M$)
we consider the alternating sum
\begin{equation}\label{eq:alt_sum}
  \sum_{\e_1,\e_2,\cdots,\e_d = 0,1}
    (-1)^{\e_1 + \e_2 + \cdots + \e_d}
    [M+\e_1 E_{i_1 i_1}+\e_2 E_{i_2 i_2}+\cdots+\e_d E_{i_d i_d}]
  \in\cS
\end{equation}
where $E_{ii}$ is the matrix of the same size as $M$
with $(i,i)$-entry 1 and the others $0$.
Let $\cS_d$ be the vector subspace spanned by the alternating sums
\eqref{eq:alt_sum} for all $M$ and $i_1,i_2,\dots,i_d$.
We call a map $f:\cS\to\Q$ a finite type invariant of degree $d$
if it vanishes on $\cS_{d+1}$.
Note that the space of all the finite type invariants of degree
$d$ is the dual space of $\cS/\cS_{d+1}$.
\par
%%%%%%%%%%%%%%%%%%%%%%%%%%%%%%%%%%%%%%%%%%%%%%%%%%%%%%%%%%%%%%%%%%%%%
\begin{rem}
It is easy to see that the transformation
$M\Leftrightarrow M+E_{i,i}$ is equivalent to the algebraic unknotting
operation introduced by the first-named author \cite{Murakami:QUEAG90}.
Therefore it follows that $\cS/\cS_1\cong\Q$ (generated by the
zero by zero matrix).
\end{rem}
%%%%%%%%%%%%%%%%%%%%%%%%%%%%%%%%%%%%%%%%%%%%%%%%%%%%%%%%%%%%%%%%%%%%%
The aim of this paper is to describe the graded space
$\oplus_{d=0}^{\infty}\cS_d/\cS_{d+1}$ in terms of web diagrams.
%%%%%%%%%%%%%%%%%%%%%%%%%%%%%%%%%%%%%%%%%%%%%%%%%%%%%%%%%%%%%%%%%%%%%
\begin{thm}\label{thm:main}
The graded vector space
$\oplus_{d=0}^{\infty} \cS_d/\cS_{d+1}$
is isomorphic to the polynomial ring $\Q[w_2, w_4, w_6, \cdots]$
as graded vector spaces, where $w_n$ is an indeterminate of degree
$n$.
\end{thm}
%%%%%%%%%%%%%%%%%%%%%%%%%%%%%%%%%%%%%%%%%%%%%%%%%%%%%%%%%%%%%%%%%%%%%
\begin{rem}
The indeterminate $w_{k}$ corresponds to the wheel $\omega_{k}$
described in the next section.
\end{rem}
%%%%%%%%%%%%%%%%%%%%%%%%%%%%%%%%%%%%%%%%%%%%%%%%%%%%%%%%%%%%%%%%%%%%%
\section{Web diagrams and clasper theory}
In this section we briefly describe web diagrams, Vassiliev
invariants, and K.~Habiro's clasper theory.
\par
A web diagram is a (possibly disconnected) uni-trivalent graph
with univalent vertices attached to the oriented $S^1$.
We use dashed lines for uni-trivalent graphs.
Let $\cA(S^1)$ be the vector space over $\Q$ spanned by all the web
diagrams modulo the AS (antisymmetry), STU, IHX, and FI
(framing independence) relations as described below.
\begin{align*}
  \text{AS relation} \qquad&\fig{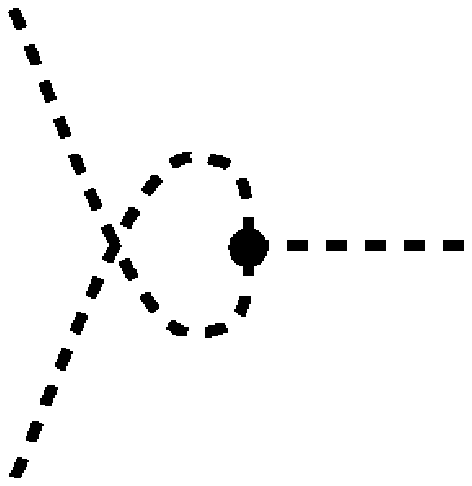}\qquad
  =\qquad-\quad\fig{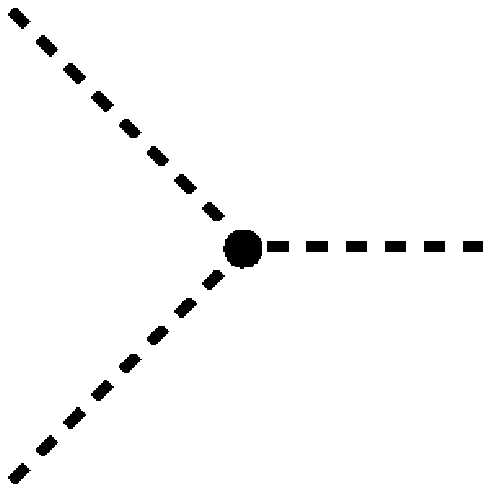}\quad ,
  \\[5mm]
  \text{IHX relation}\qquad&\fig{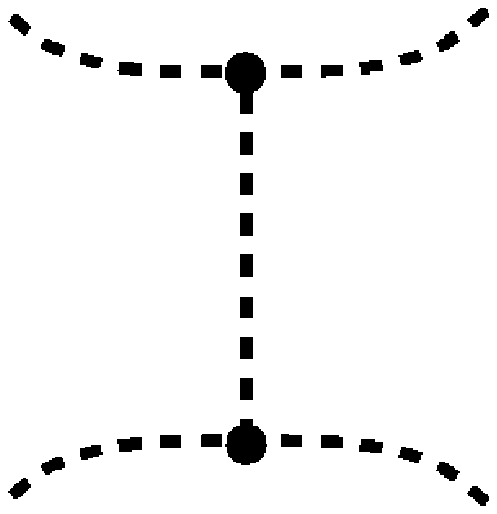}\qquad
  =\qquad\fig{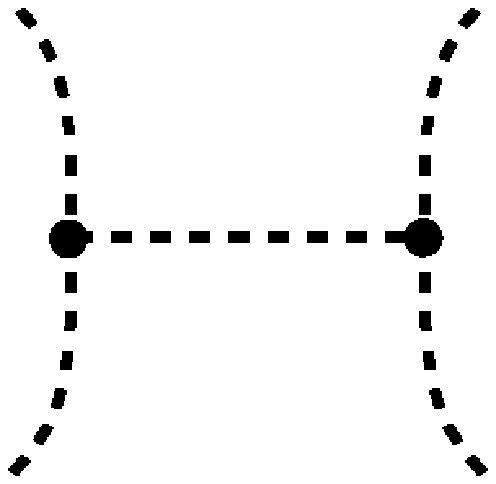}\quad-\quad\fig{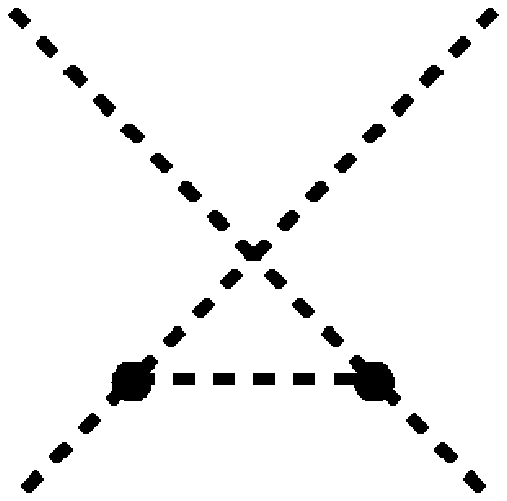}\quad ,
  \\[5mm]
  \text{STU relation}\qquad&\fig{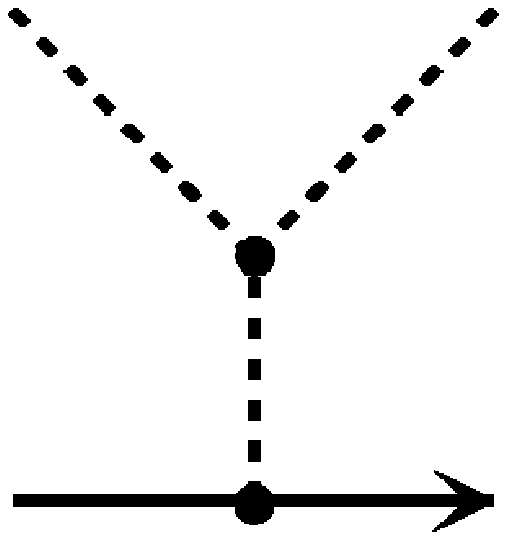}\qquad
  =\qquad\fig{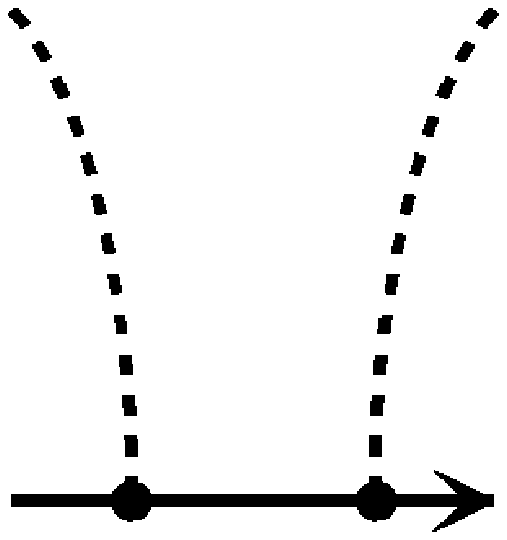}\quad-\quad\fig{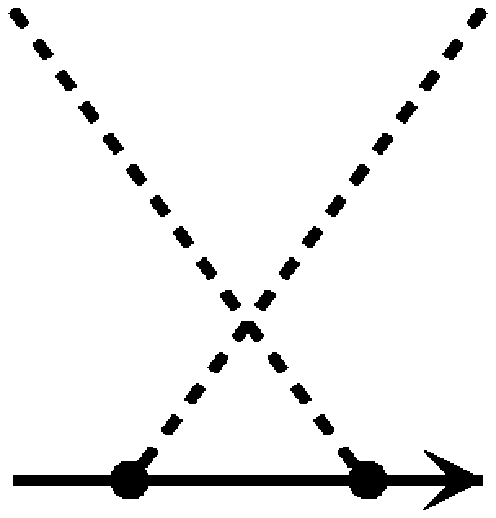}\quad ,
  \\[5mm]
  \text{FI relation} \qquad&\fig{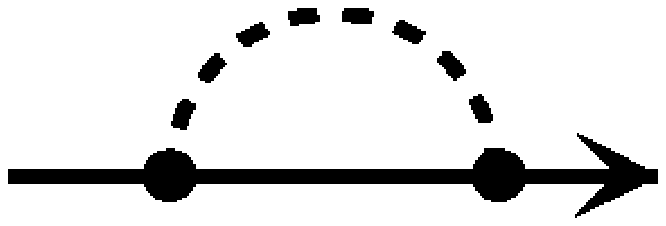}\qquad=\qquad 0.
\end{align*}
Note that all the vertices where three dashed edges meet are oriented
counter-clockwise and in the AS relation we respect this orientation.
Here we add the FI relation since we are studying unframed knots.
\par
The space $\cA(S^1)$ has an algebra structure with respect to
the connected sum of $S^1$.
It is known to be generated as an algebra by web diagrams with
connected dashed uni-trivalent graphs.
Some of the generators are known as ``wheels''
(see
\cite{Chmutov/Varchenko:TOPOL97,
      Kricker:JKNOT97,
      Kricker/Spence/Aitchison:JKNOT97,
      Bar-Natan/Garoufalidis/Rozansky/Thurston:1997})
denoted by $\omega_{2n}$;
see Figure~1 for their definitions.
(It is denoted by $\tau_{\{2n\}}$ in
\cite{Kricker:JKNOT97,Kricker/Spence/Aitchison:JKNOT97}.
Note that the definition of wheels in
\cite{Bar-Natan/Garoufalidis/Rozansky/Thurston:1997}
is different from ours.)
\begin{gather*}
  \begin{matrix}\fig{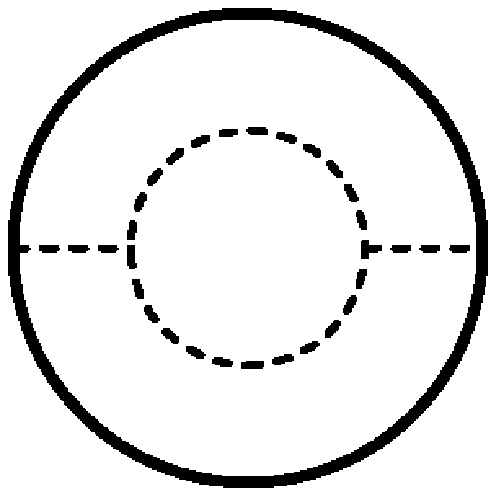} \\ \omega_{2}\end{matrix}\quad,\qquad
  \begin{matrix}\fig{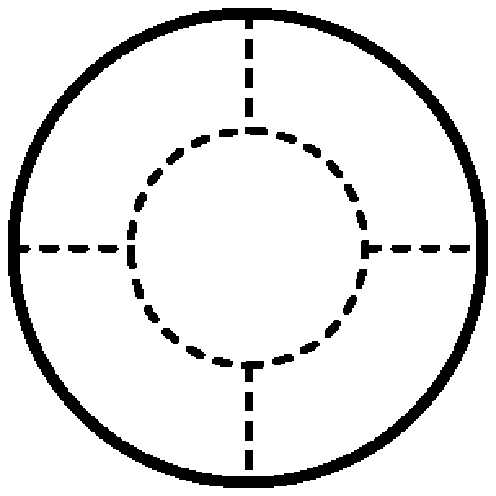} \\ \omega_{4}\end{matrix}\quad,\qquad
  \begin{matrix}\fig{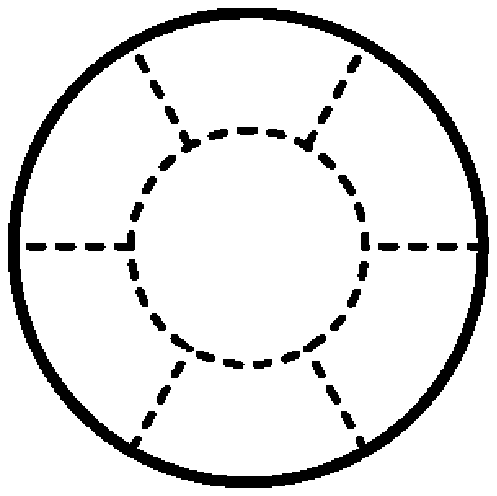} \\ \omega_{6}\end{matrix}\quad,\dots .
  \\
  \text{Figure~1.}
\end{gather*}
The other generators are web diagrams whose dashed parts have negative
Euler characteristics.
\par
A map from $\cA(S^1)$ to $\Q$ is called a weight system.
One can construct a weight system from coefficients
of the Alexander-Conway polynomial \cite{BarNatan/Garoufalidis:INVEM96}.
Moreover we can characterize a weight system which comes from
(a sum of products of) coefficients of the Alexander-Conway
polynomial as follows due to A.~Kricker, B.~Spence and I.~Aitchison
\cite{Kricker/Spence/Aitchison:JKNOT97} (see also
\cite[Lemma 2.11]{Kricker:JKNOT97}).
\begin{lem}\label{lem:KSA}
Let $W$ be a weight system vanishing on web diagrams with a
connected dashed component of negative Euler characteristic.
Then it equals a sum of products of weight systems coming from
coefficients of the Alexander-Conway polynomial.
\end{lem}
\par
Next we review the definition of Vassiliev invariants of knots;
see \cite{BarNatan:TOPOL95} for detailed definition for example.
Let $\cK$ be the vector space over $\Q$ spanned by the isotopy classes
of knots in $S^3$.
Let $\cK_d$ be the vector subspace of $\cK$ spanned by singular knots
with $d$ double points;
a singular knot can be regarded as an element of $\cK$
in the sense that a double point is regarded as the difference of
the positive and negative crossing as follows.
\begin{equation*}
  \fig{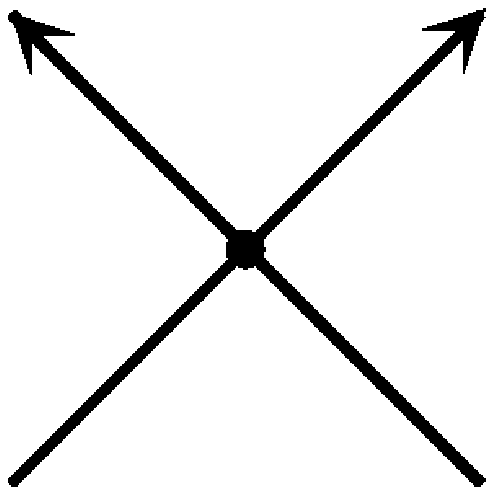}\qquad=\qquad\fig{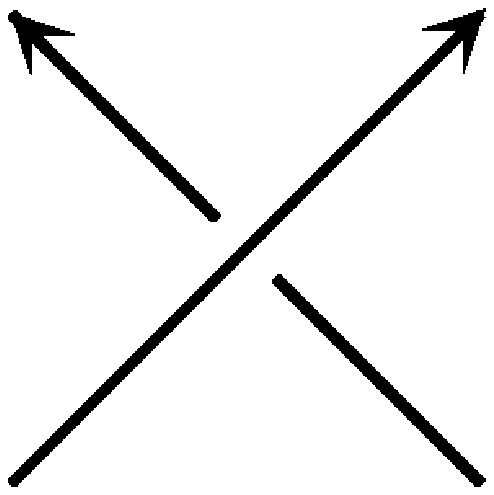}\qquad-\qquad\fig{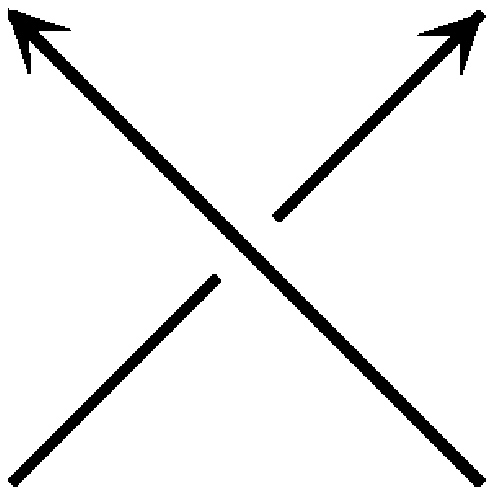}\quad .
\end{equation*}
A linear map $v : \cK \to \Q$ is called a {\it Vassiliev invariant}
(or a {\it finite type invariant}) of degree $d$
if the map $v$ vanishes in $\cK_{d+1}$.
\par
From the STU relation any web diagram can be represented as
a linear combination of chord diagrams (web diagrams without
dashed trivalent vertices).
A map $\varphi:\cA(S^1)^{(d)}\to\cS_{d}/\cS_{d+1}$ is defined as
follows, where $\cA(S^1)^{(d)}$ is the degree $d$ part
of $\cA(S^1)$ (the degree of a web diagram is half the number of
vertices).
For a chord diagram $E$, we define $\varphi(E)$ to be the singular
knot whose double points correspond to chords (dashed lines).
This is well-defined modulo $\cS_{d+1}$.
Now we extend $\varphi(D)$ linearly to the vector space spanned
by chord diagrams.
Finally we extend it to $\cA(S^1)^{(d)}$ by using the STU relation.
It is not hard to check that $\varphi$ is well-defined as a map
to $\cS_{d}/\cS_{d+1}$.
Moreover it is known that $\varphi$ is an isomorphism by Kontsevich's
integral (see \cite{BarNatan:TOPOL95}).
\par
We have a natural linear map $s: \cK \to \cS$ which takes a knot to
the S-equivalence class of a Seifert matrix for the knot.
Given a singular knot, we express it as a linear combination of knots.
We may assume that these knots differ only near double points of the
singular knot.
Since we can choose Seifert surfaces for these knots in such a way
that their Seifert matrices are as in \eqref{eq:alt_sum},
the image of $\cK_d$ by $s$ is in $\cS_d$.
Hence the map $s$ induces the map $\cK_d/\cK_{d+1}\to\cS_d/\cS_{d+1}$,
which we also denote by $s$.
\par
From Lemma~\ref{lem:KSA}, we have the following corollary to
Theorem~\ref{thm:main}.
%%%%%%%%%%%%%%%%%%%%%%%%%%%%%%%%%%%%%%%%%%%%%%%%%%%%%%%%%%%%%%%%%%%%%
\begin{cor}\label{cor}
Any Vassiliev invariant of knots which factors a finite type invariant
of Seifert matrices is equal to a linear sum of coefficients of
the Alexander polynomial.
\end{cor}
%%%%%%%%%%%%%%%%%%%%%%%%%%%%%%%%%%%%%%%%%%%%%%%%%%%%%%%%%%%%%%%%%%%%%
\par
Now we briefly review K.~Habiro's clasper theory.
(For more detail and for more general theory, we refer the reader to
\cite{Habiro:Oiwake97}.)
\par
Let $K\in S^3$ be a knot in the $3$-space.
A clasper $G$ for $K$ is a framed uni-trivalent graph embedded in
$S^3$ with its univalent vertices on $K$ such that its interior does
not touch $K$.
So it is an embedding of a web diagram in $S^3$ with framing.
We use the blackboard framing to describe claspers.
The degree of the clasper $G$ is half the number of vertices and
denoted by $\deg(G)$.
\begin{gather*}
  \begin{matrix}\fig{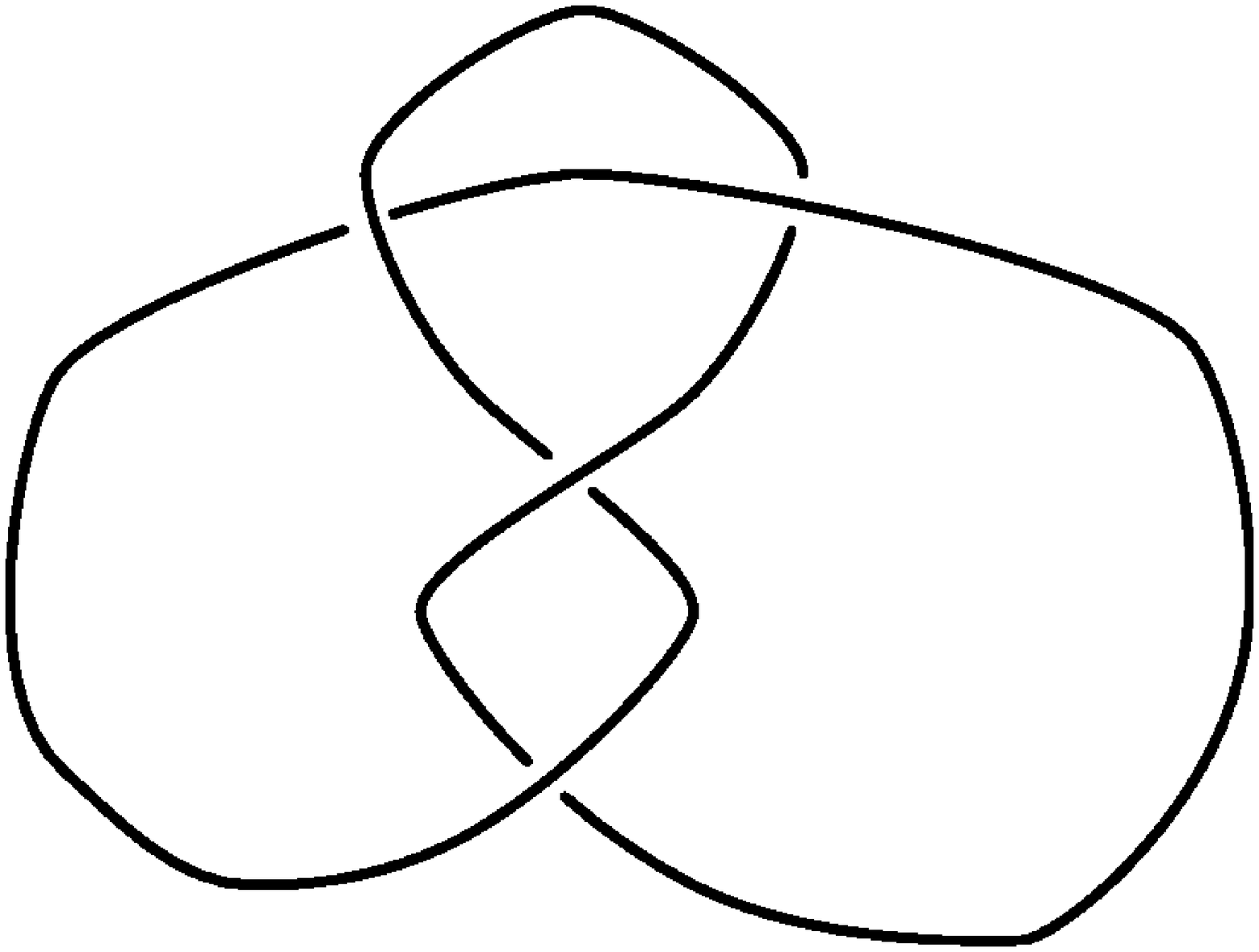}\\[10mm] \text{a knot $K$}\end{matrix}
  \qquad
  \begin{matrix}\fig{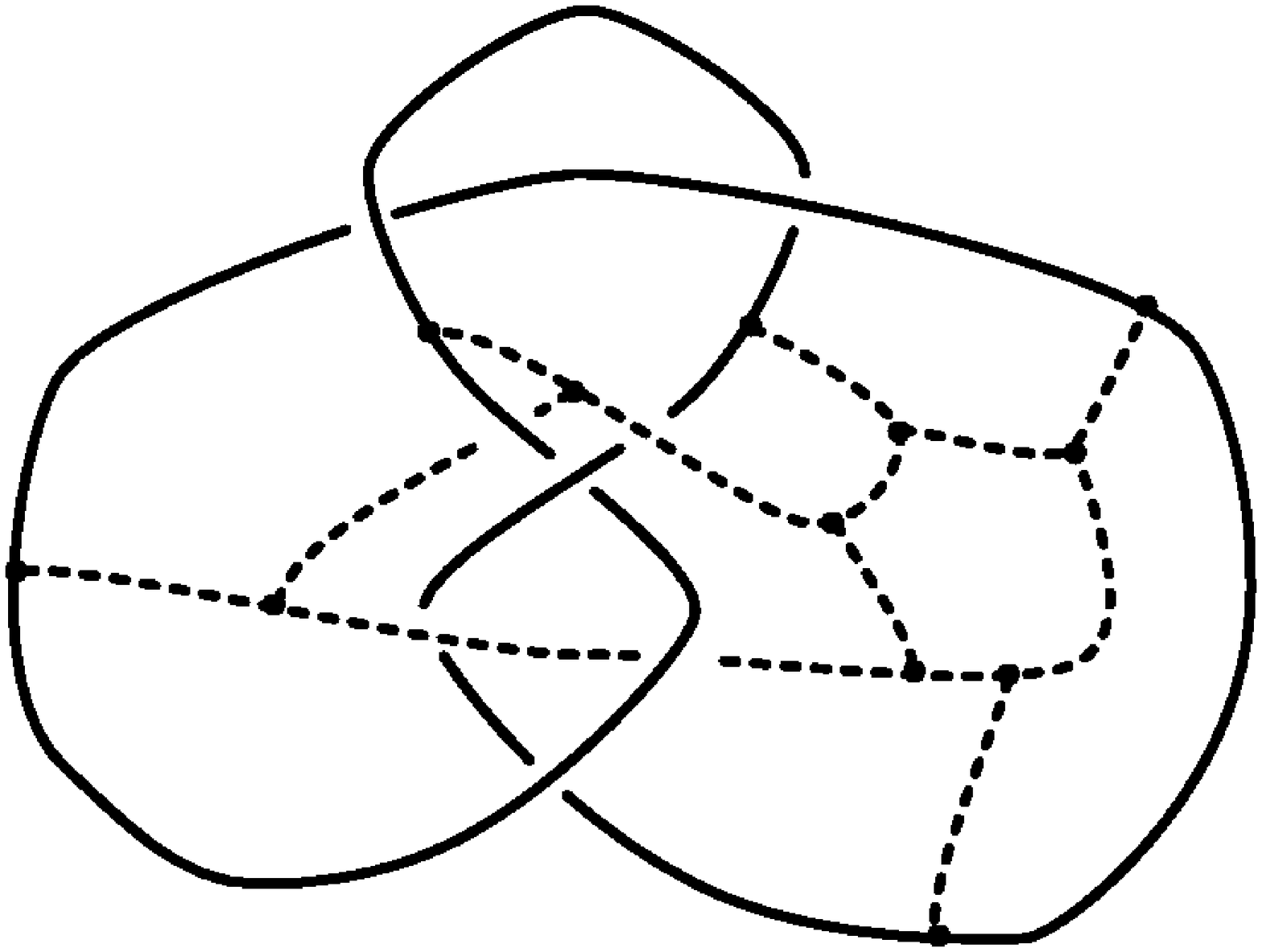}\\[10mm] \text{a clasper for $K$}\end{matrix}
  \\
  \text{Figure~2.}
\end{gather*}
\par
By $K_{G}$ we denote the knot obtained as follows.
First we replace each trivalent vertex with Borromean rings and
each univalent vertex with a clasp (see Figure~3).
\begin{gather*}
  \raisebox{-5mm}{\fig{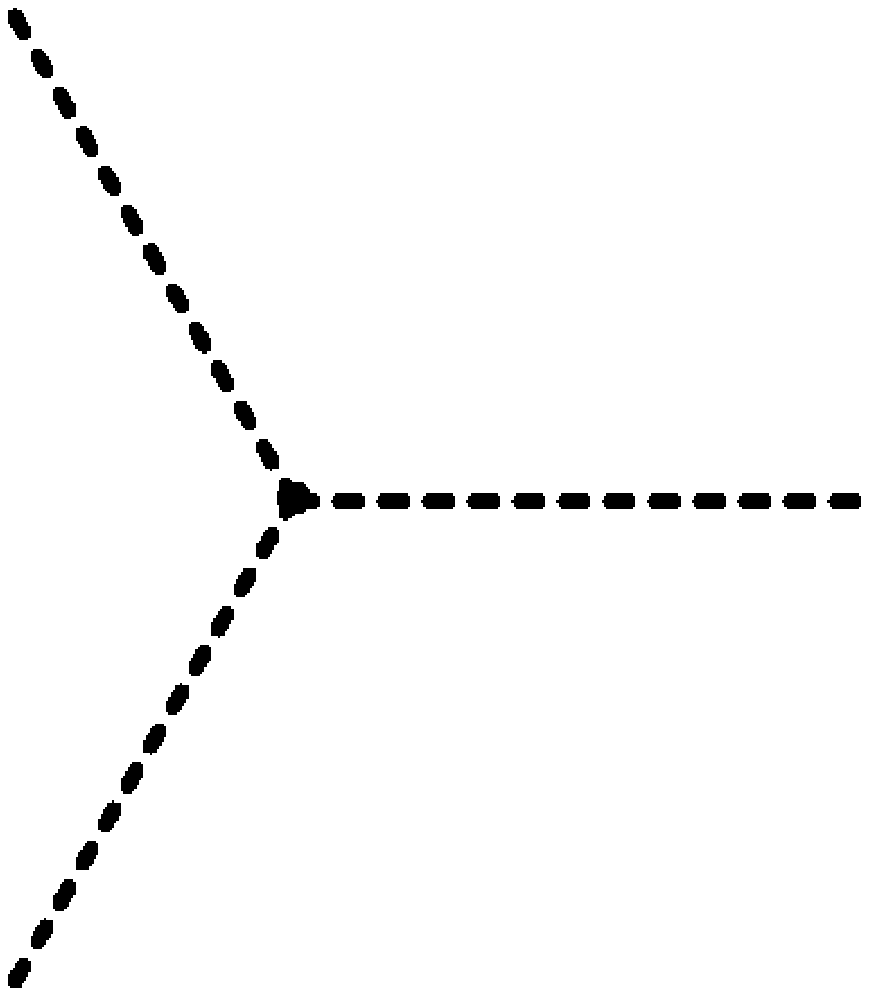}}\qquad\Rightarrow\qquad
  \raisebox{-5mm}{\fig{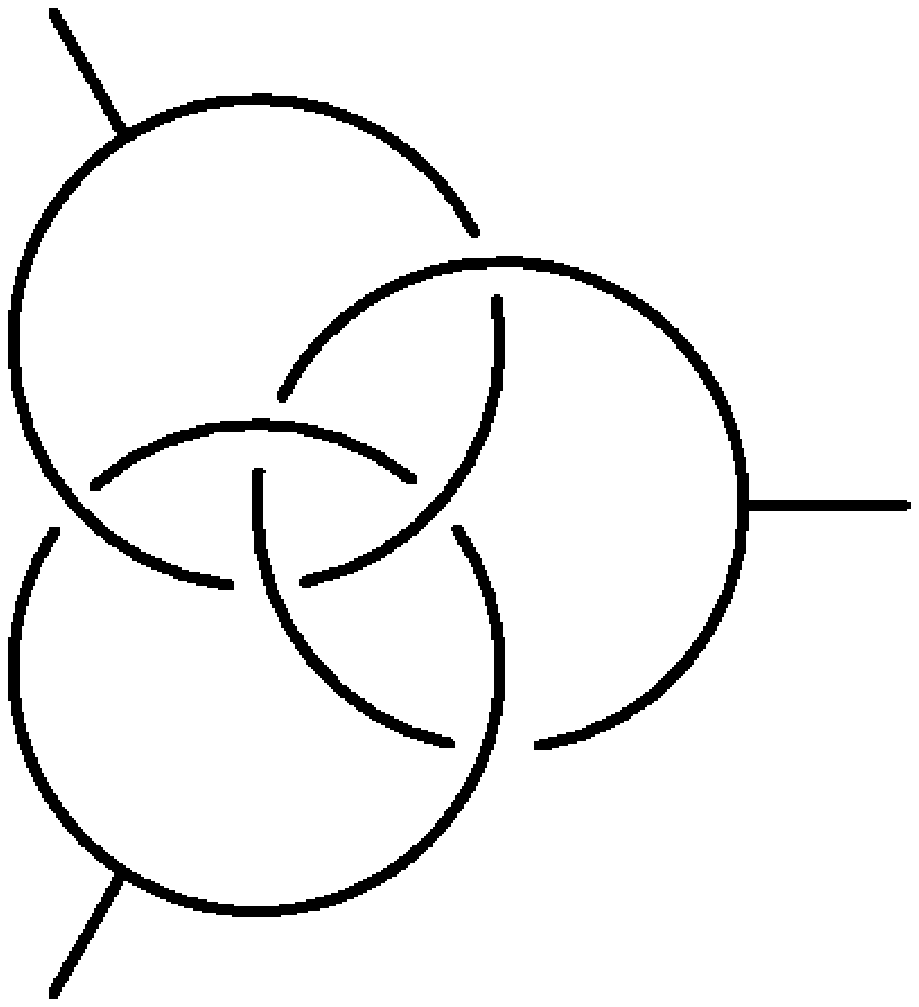}}
  \notag
  \\
  \raisebox{-3mm}{\fig{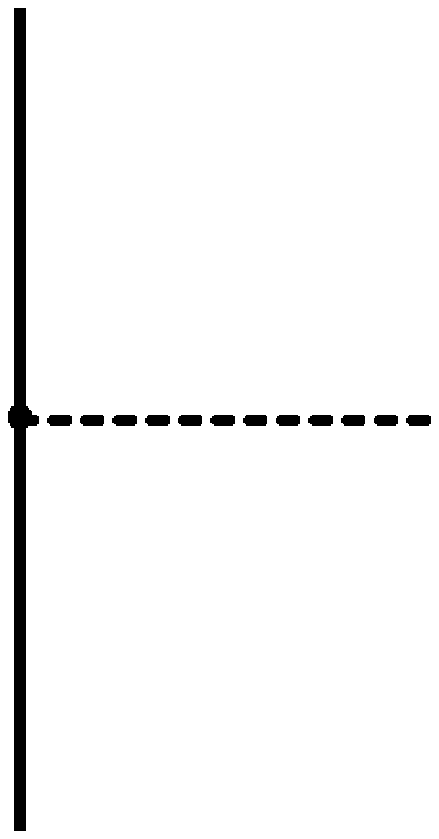}}\qquad\Rightarrow\qquad
  \raisebox{-3mm}{\fig{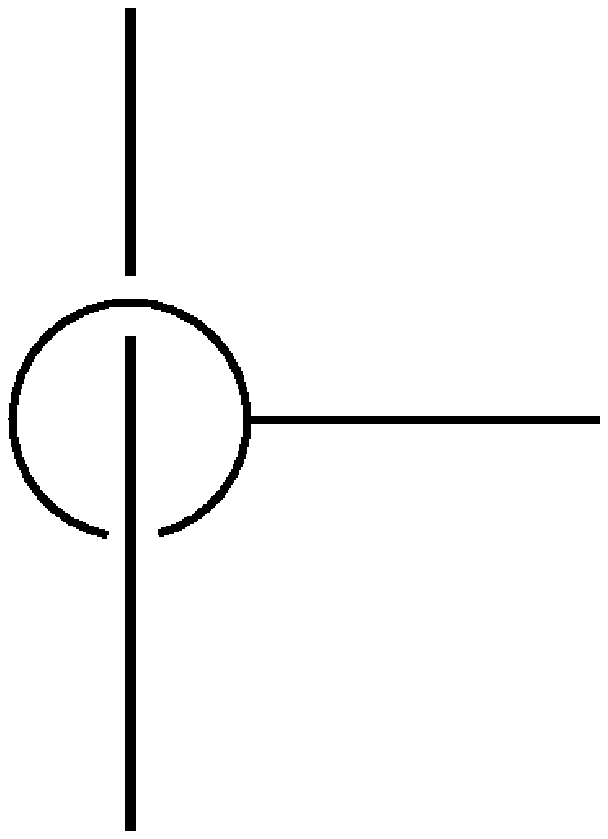}}
  \\
  \text{Figure~3.}
\end{gather*}
Next we replace each edge of the resulting trivalent graph with a Hopf
link as in Figure~4.
\begin{gather*}
  \raisebox{-6mm}{\fig{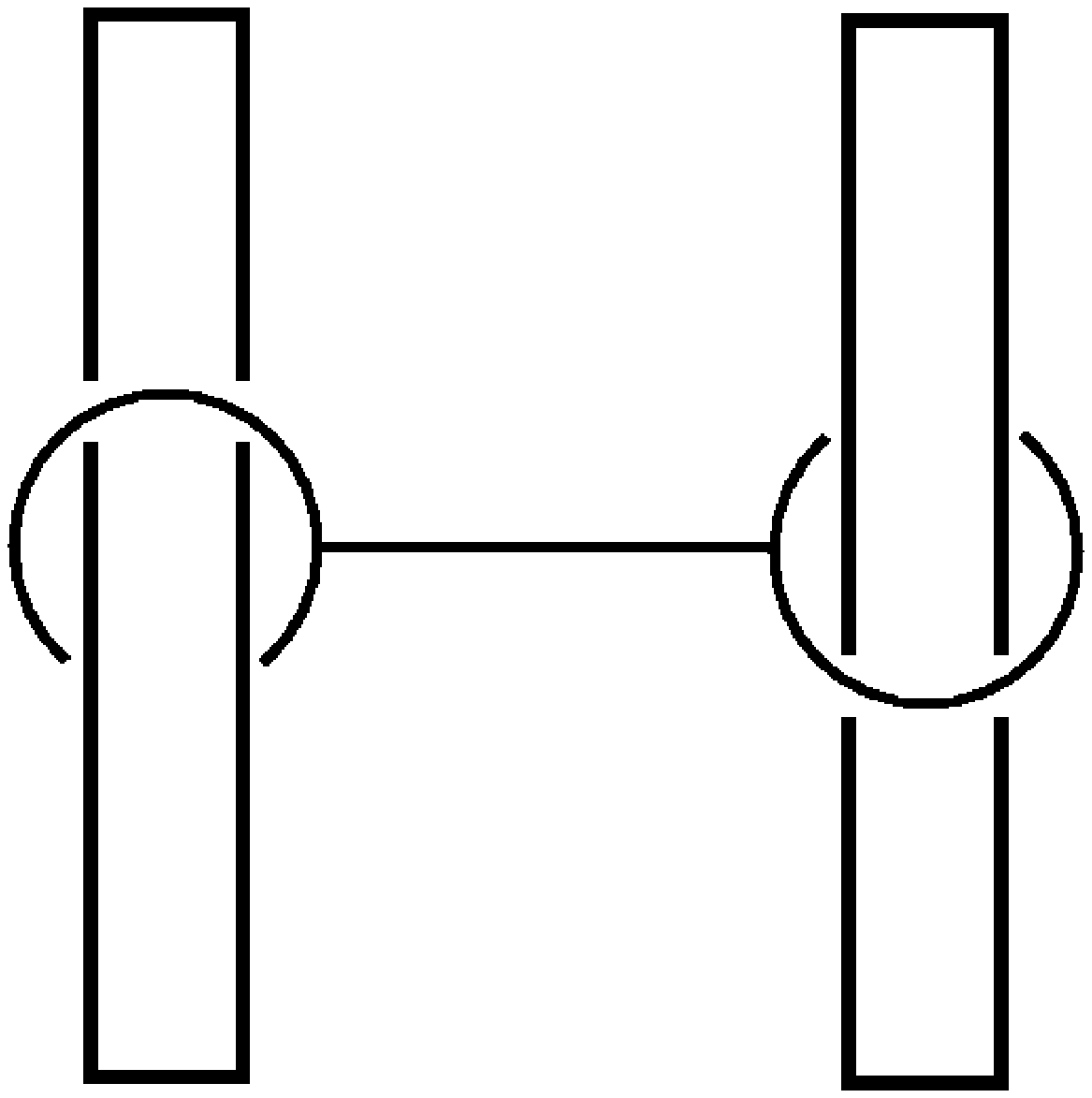}}\qquad\Rightarrow\qquad
  \raisebox{-6mm}{\fig{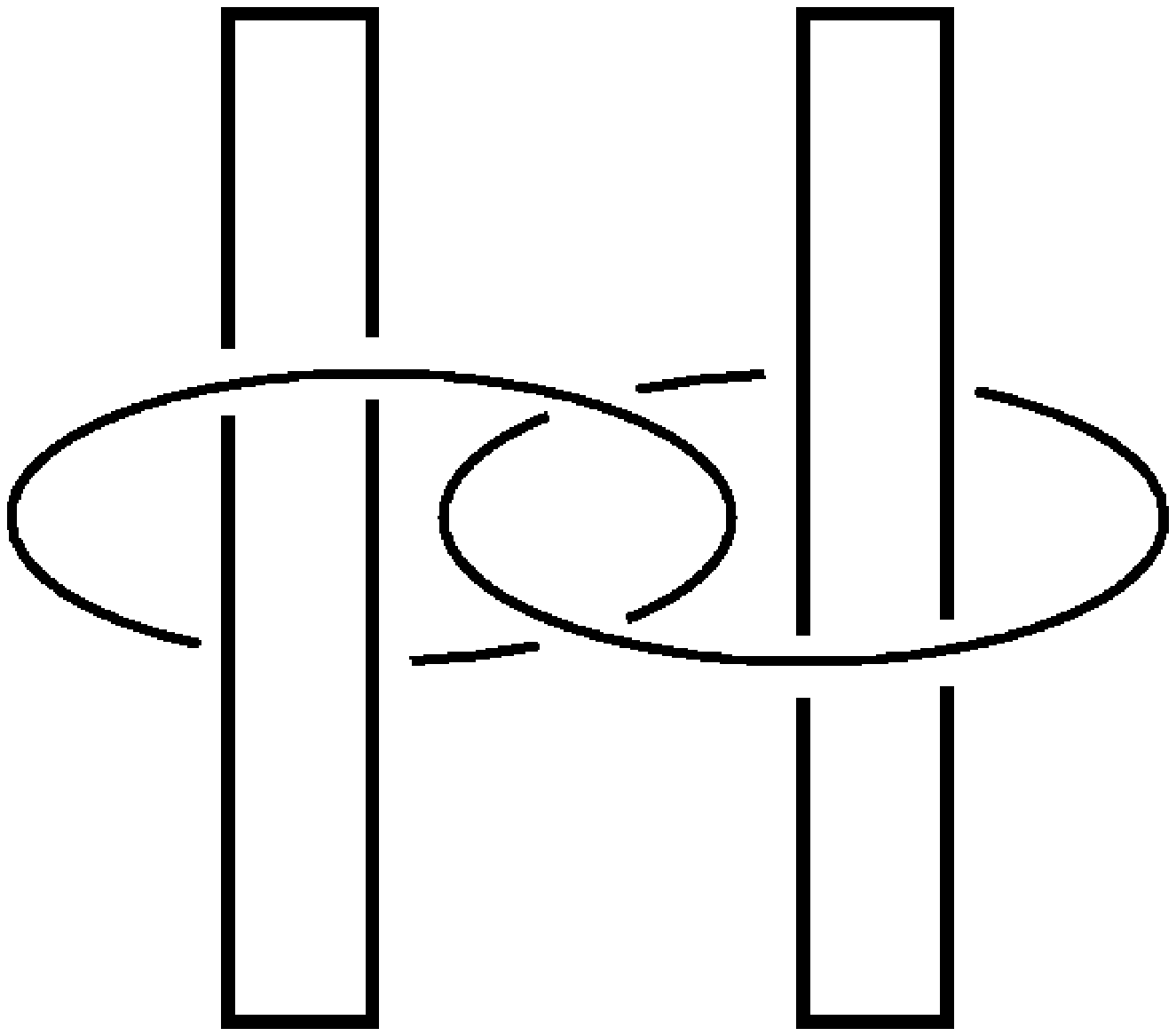}}
  \\
  \text{Figure~4.}
\end{gather*}
We denote the resulting link by $L_{G}$.
Let $\cL_{G}$ be the framed link obtained from $L_{G}$ with
every framing $0$.
Then it is easy to see that $S^3$ surgered along $\cL_{G}$ is again
$S^3$.
Now $K_{G}$ is defined to be the knot in this surgered $S^3$.
\par
For a pair $(K,G)$ of a knot and a clasper for it, we put
\begin{equation*}
  e(K,G)=
  \sum_{\emptyset\subseteq G'\subseteq G}(-1)^{\deg(G)-\deg(G')}K_{G}
  \in\cS,
\end{equation*}
where $G'$ runs over all the unions of connected components in $G$
including the empty set.
Let $\cG_d$ be the vector space spanned by the pairs $(K,G)$
with $\deg(G)=d$.
Then $e$ defines a map from $\cG_d$ to $\cK_d$, which we also
denote by $e$.
One of the main result of Habiro's clasper theory is the following.
%%%%%%%%%%%%%%%%%%%%%%%%%%%%%%%%%%%%%%%%%%%%%%%%%%%%%%%%%%%%%%%%%%%%%
\begin{thm}[\cite{Habiro:GEOTO00,Habiro:Oiwake97}]\label{thm:surjective}
Let $\gamma:\cG_{d}\to\cA(S^1)^{(d)}$ be the map forgetting
the embedding.
Then the following diagram commutes.
\begin{equation*}
\begin{CD}
  \cG_{d}        @>{e}>>              \cK_{d}               \\
  @V{\gamma}VV                        @VV\text{projection}V \\
  \cA(S^1)^{(d)} @>{\cong}>{\varphi}> \cK_{d}/\cK_{d+1}
\end{CD}
\end{equation*}
Moreover  $e:\cG_d\to\cK_d$ is surjective.
\end{thm}
\par
We note that M.N.~Goussarov has obtained a similar result using
``Y-graphs'' \cite{Goussarov:1998}.
%%%%%%%%%%%%%%%%%%%%%%%%%%%%%%%%%%%%%%%%%%%%%%%%%%%%%%%%%%%%%%%%%%%%%
\section{Proof of Theorem~\ref{thm:main}}
%%%%%%%%%%%%%%%%%%%%%%%%%%%%%%%%%%%%%%%%%%%%%%%%%%%%%%%%%%%%%%%%%%%%%
In this section we prove Theorem~\ref{thm:main}.
\par
From the previous section we have the following maps
\begin{equation*}
  \Q[\omega_2,\omega_4,\omega_6,\cdots]^{(d)}\xrightarrow{i}
  \cA(S^1)^{(d)}\xrightarrow{\varphi}
  \cK_d/\cK_{d+1}\xrightarrow{s}
  \cS_d/\cS_{d+1}.
\end{equation*}
Here $\Q[\omega_2,\omega_4,\omega_6,\cdots]^{(d)}$ is the degree
$d$ part of $\Q[\omega_2,\omega_4,\omega_6,\cdots]$.
In the following we will show that the composed map
$s\circ\varphi\circ i$ is an isomorphism.
\par
First we show that $s\circ\varphi\circ i$ is surjective.
\begin{proof}[Surjectivity]
Since $s$ is surjective by definition and $\varphi$ is an isomorphism,
for any $[M]\in\cS_{d}/\cS_{d+1}$ there exists
$x\in\cA(S^1)^{(d)}$ with
$s\circ\varphi(x)=[M]$.
As we described before $\cA(S^1)^{(d)}$ is spanned by connected sums
of web diagrams whose dashed part are wheels with even spokes or
uni-trivalent graphs with negative Euler characteristics.
Therefore it is sufficient to show that a web diagram with a dashed
part of negative Euler characteristic is taken to $0$ by the map
$s\circ\varphi$.
\par
Let $D$ be a web diagram with a dashed part of negative Euler
characteristic.
Then $\varphi(D)$ is expressed as the difference between
a knot $K$ and a knot $K'$ which is obtained from $K$
by surgery along a set of claspers corresponding to $D$ from
Theorem~\ref{thm:surjective}.
We easily see that near a univalent vertex $K'$ is obtained by
$0$-surgery along the links shown in Figure~5.
\begin{gather*}
  \raisebox{-5mm}{\figm{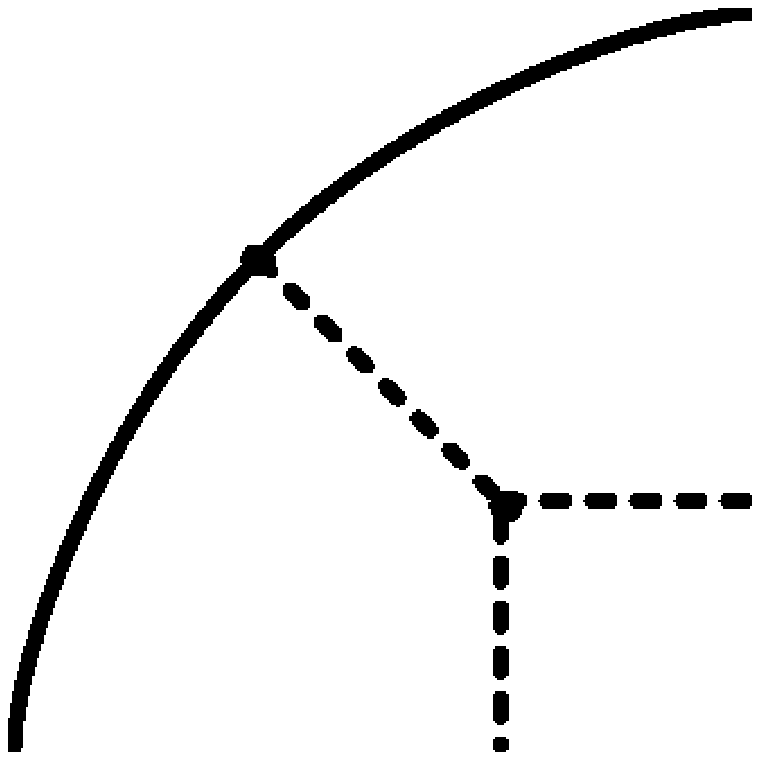}}\qquad\Rightarrow\qquad
  \raisebox{-5mm}{\figm{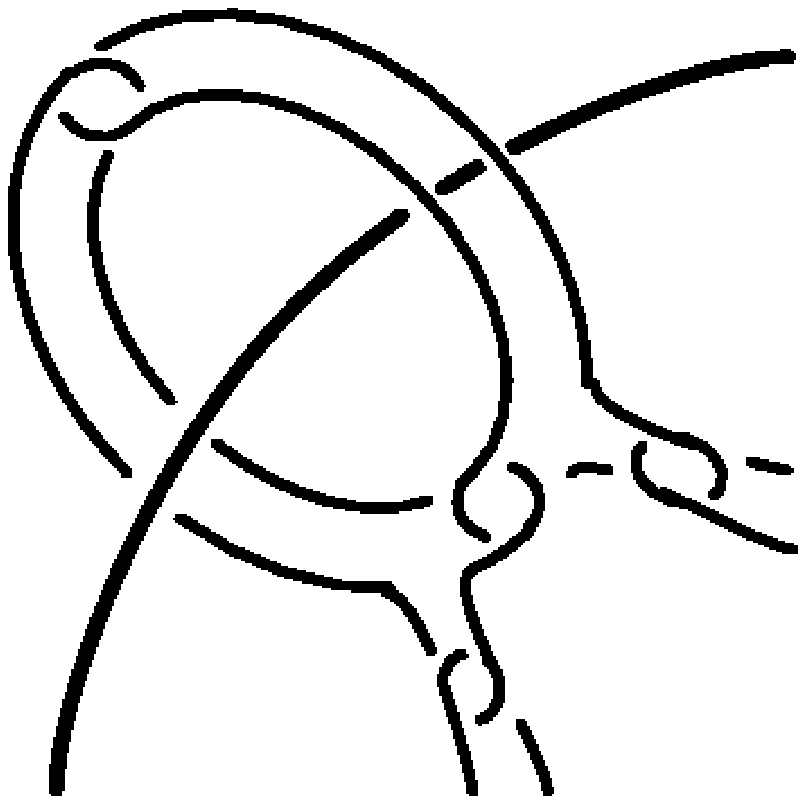}}
  \\
  \text{Figure~5.}
\end{gather*}
Hence Seifert surfaces of $K$ and $K'$ can be constructed
as shown in Figure~6.
\begin{gather*}
  \begin{matrix}\figl{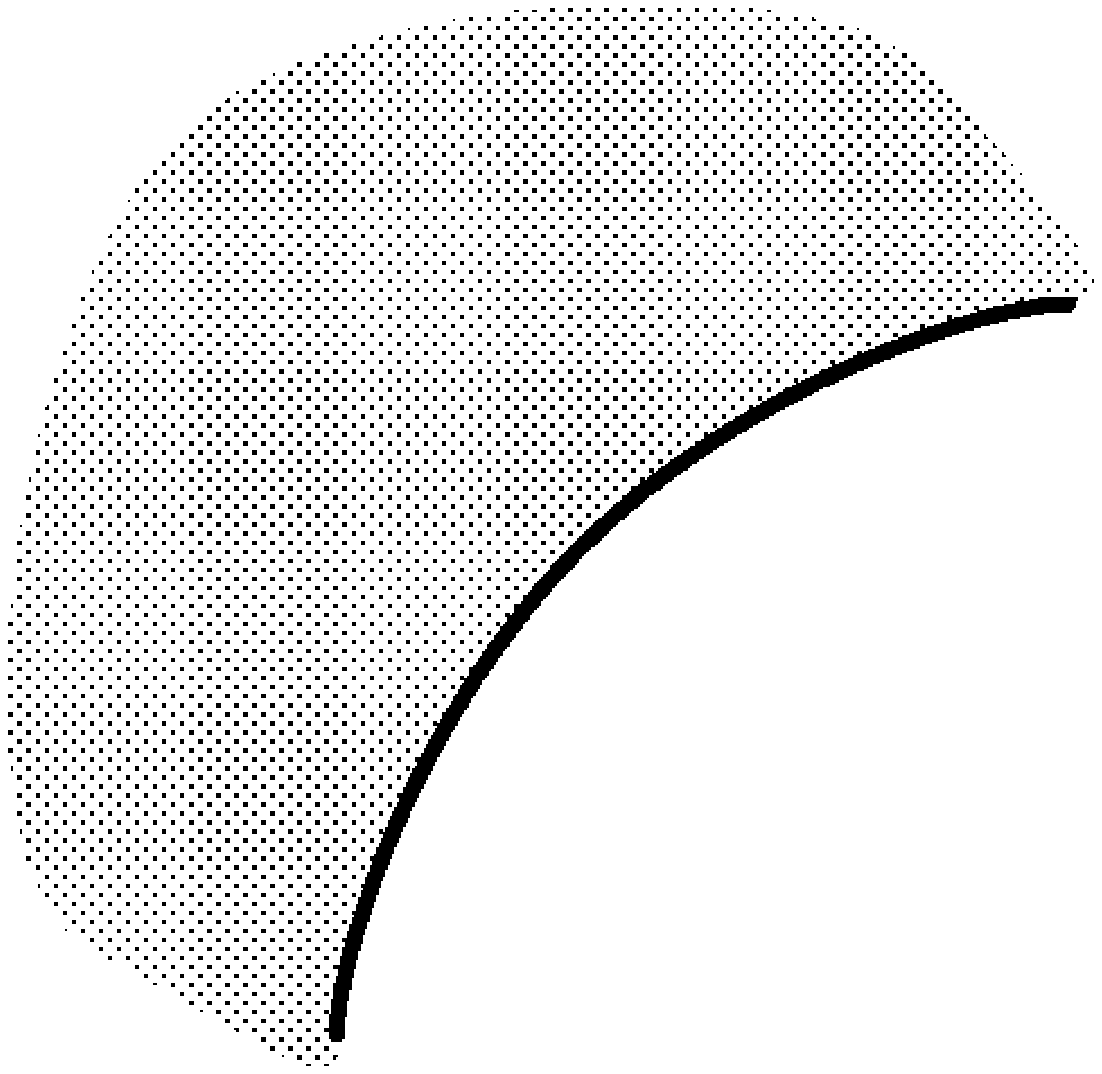}\\[10mm] K\end{matrix}
  \qquad
  \begin{matrix}\figl{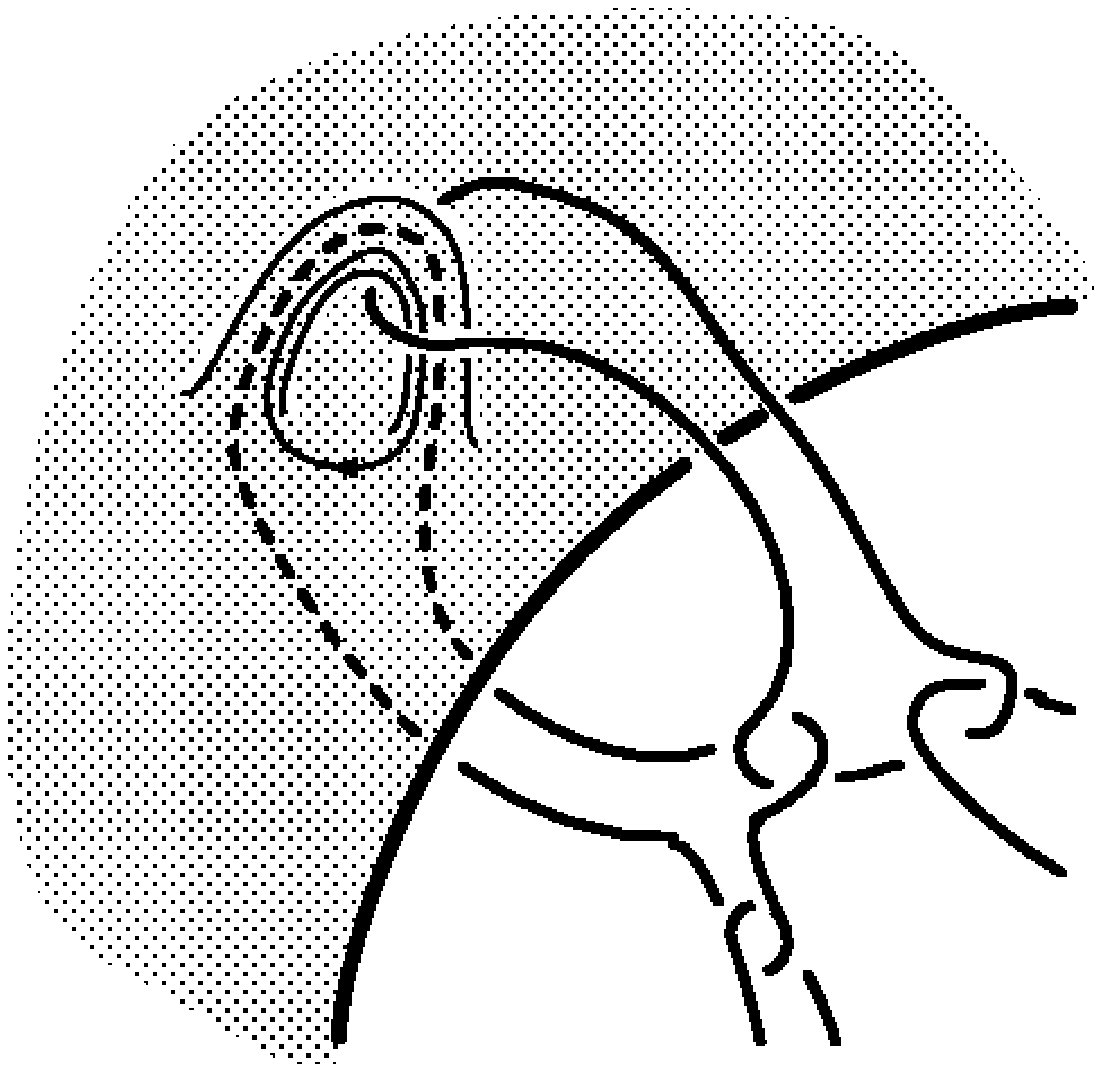}\\[10mm] K'\end{matrix}
  \\
  \text{Figure~6. Seifert surfaces for $K$ and $K'$.}
\end{gather*}
\par
Now since $D$ has a dashed uni-trivalent graph with negative
Euler characteristic, there exists a dashed trivalent vertex
which is not next to a univalent vertex on $S^1$.
Therefore the one-handle in the Seifert surface of $K'$ shown in
Figure~6 is homologically trivial.
(Note that the homology generator indicated as the arrowed circle in Figure~6
is null-homologous in such a case.
Figure~7 shows a 2-chain bounding the generator.)
\begin{gather*}
  \figl{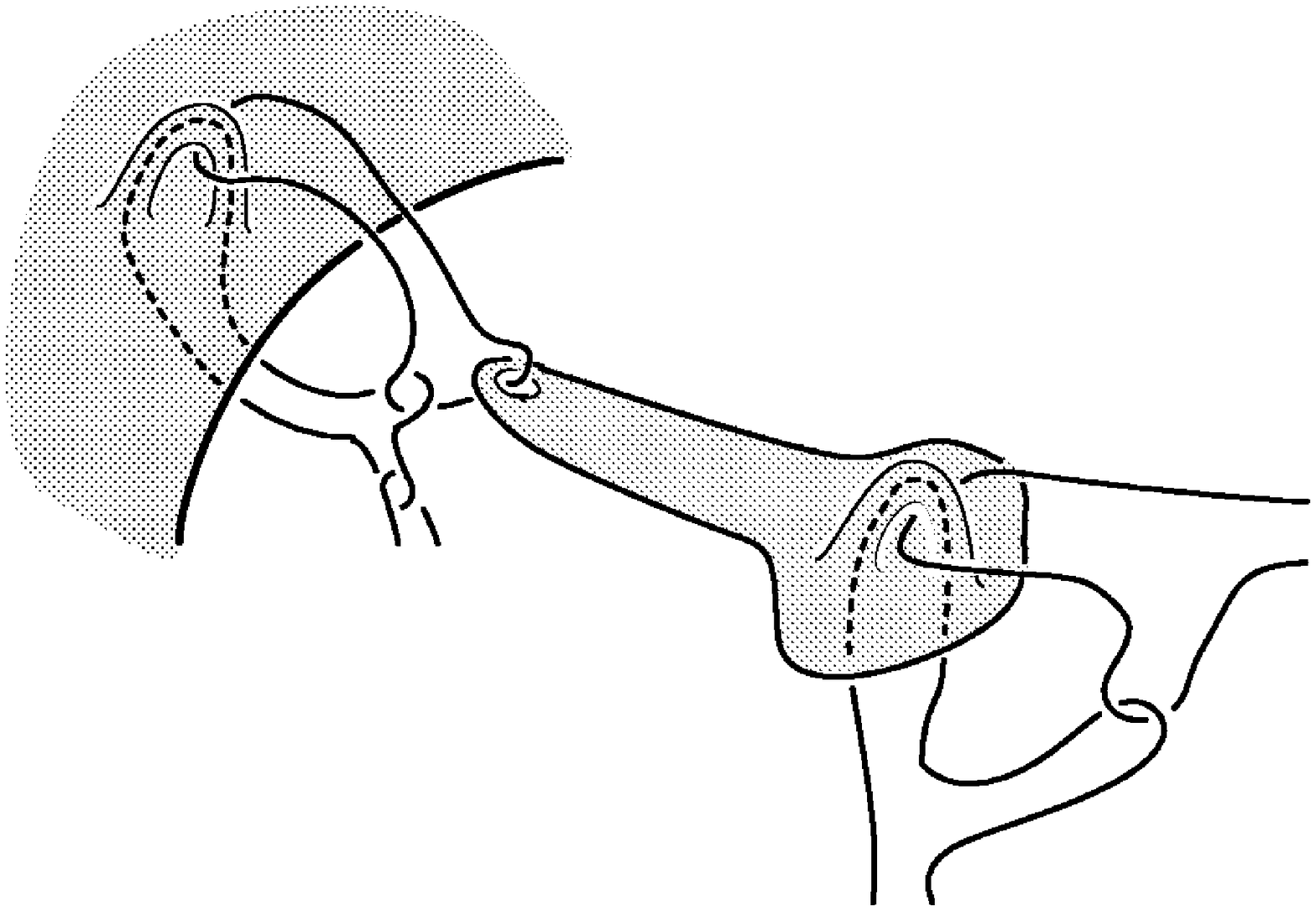}
  \\
  \text{Figure~7.}
\end{gather*}
Hence the Seifert matrices of $K$ and $K'$ are S-equivalent.
This implies that the image of $D$ by the map $s \circ \varphi$
vanishes.
\end{proof}
Next we show that $s\circ\varphi\circ i$ is injective.
\begin{proof}[Injectivity]
Suppose for a contradiction that there exists a non-zero web diagram
$D\in\Q[\omega_2,\omega_4,\omega_6\dots]^{(d)}$
we have $s\circ\varphi\circ i(D)=0$.
Now we define a weight system $W:\cA(S^1)^{(d)}\to\Q$ by
$W(D)=1$ and $W(E)=0$ if $E\ne D$ (and extend it linearly to
$\cA(S^1)^{(d)}$).
This is well-defined since the set of wheels (with even spokes)
are linearly independent with respect to the AS, IHX, STU and FI
relations (see the remark after Definition~4.8 in
\cite{Kricker/Spence/Aitchison:JKNOT97}).
Now from Lemma~\ref{lem:KSA} the weight system $W$ is derived from
coefficients of the Alexander-Conway polynomial.
\par
Since the Alexander-Conway polynomial can be obtained via
Seifert matrices, $W$ factors $S_d/S_{d+1}$ and so
$W\circ i(D)=0$.
But this is a contradiction since $W(i(D))=1$ by the definition
of $W$.
\end{proof}
\bibliography{mrabbrev,hitoshi}
\bibliographystyle{amsplain}
\end{document}